\documentclass[12pt]{amsart}
\usepackage{amscd,amssymb}
\usepackage[arrow,matrix]{xy}
\usepackage[plainpages,backref,urlcolor=blue]{hyperref}
\usepackage{mathdots}
\usepackage{mathrsfs}

\theoremstyle{plain}
\newtheorem{thm}[subsection]{Theorem}

\newtheorem{cor}[subsection]{Corollary}

\theoremstyle{definition}

\newtheorem{definition}[subsection]{Definition}
\newtheorem{ex}[subsection]{Example}

\numberwithin{equation}{section}
\newcommand{\C}{\mathbb{C}}

\newcommand{\PP}{\mathbb{P}}

\begin{document}
\date{}

\title [On homogeneous polynomials  determined by their Jacobian ideal]
{On homogeneous polynomials  determined by their Jacobian ideal}

\author[ZHENJIAN WANG]{ ZHENJIAN WANG  }
\address{Univ. Nice Sophia Antipolis, CNRS,  LJAD, UMR 7351, 06100 Nice, France.}
\email{wzhj@mail.ustc.edu.cn}



\begin{abstract} We investigate which homogeneous polynomials are determined by their Jacobian ideals, and extend and complete previous results due to J. Carlson and Ph. Griffiths, K. Ueda and M. Yoshinaga, and A. Dimca and E. Sernesi.

\end{abstract}

\maketitle


\section{Introduction} \label{sec:intro}

Let $S=\mathbb{C}[x_0,\cdots,x_n]$ be the graded ring of polynomials in $(n+1)$ variables with complex coefficients. For a homogeneous polynomial $f \in S$ of degree $d$, we denote by $J_f$ the homogeneous ideal in $S$ generated by the first order partial derivatives
$$f_i=\frac{\partial f}{\partial x_i},$$
and by $E(f)$  the homogeneous component $J_{f,d-1}$, which is just the vector space spanned by $f_0,\cdots,f_n$.

A celebrated result due to J. Carlson and Ph. Griffiths, see section 4, (b) in \cite{CG} tells us that  a {\it generic} polynomial $f$ can be reconstructed (up to a multiplicative constant factor) from its Jacobian ideal $J_f$, or, equivalently, from the vector space $E(f)$. In this paper they conjecture that $f$ is determined by $E(f)$ when the corresponding projective hypersurface $V(f)$ ``does not possess singularities which are large in relation to the degree".

If $V(f)$  is smooth, then this conjecture was settled by K. Ueda and M. Yoshinaga in \cite{UY}. They showed that in such a case, $f$ is determined by $E(f)$ unless $f$ is of Sebastiani-Thom type, see Definition \ref{defST} below.
In a recent preprint \cite{DS}, A. Dimca and E. Sernesi have extended the result by K. Ueda and M. Yoshinaga to cover some (not very) singular plane curves, i.e. the case $n=2$.

In this note we prove the following result, confirming the above Carlson-Griffiths conjecture in the general case.

\begin{thm}\label{thm1}
Suppose given two homogeneous polynomials $f,g\in \mathbb{C}[x_0,\cdots,x_n]$ such that $V(f) \ne V(g)$, i.e. $ f\notin\mathbb{C}^*g$, of degree $d\geq3$.

If $J_f=J_g$, then either the hypersurface $V(f)\subset\mathbb{P}^n$ has a singularity of multiplicity $d-1$, or $f$ is of Sebastiani-Thom type (ST type).
\end{thm}

In particular we get the following corollaries.

\begin{cor}
If $J_f=J_g$,  $V(f) \ne V(g)$  and $V(f)$ does not admit any singularities of multiplicity $d-1$, then $f$ is of ST type.
\end{cor}

\begin{cor}
If $J_f=J_g$, $V(f) \ne V(g)$ and $V(f)$  is smooth, then $f$ is of ST type.
\end{cor}

For any $n\geq 1$ and $d \geq 3$, we show in Example \ref{existence} the existence of a pair $(f,g), f,g\in\mathbb{C}[x_0,\cdots,x_n]_d$ with $V(f) \ne V(g)$, $J_f=J_g$, and $f$  not of ST type.
For example, for $n=2$, $d=3$ we can take
$g(x,y,z)=x^2z+xy^2,
f(x,y,z)=\lambda g+x^2y,\ \ \lambda\neq0.$

In Corollary \ref{corA} we give a precise meaning of the word ``generic" in the above mentioned result by J. Carlson and Ph. Griffiths.
Both Example \ref{existence} and Corollary \ref{corA} are based on a detailed study of a special case done in subsection (4.2), which is not strictly speaking necessary for the proof of Theorem
\ref{thm1}.

Let $E'(f)=J_{f,d}$ and recall R. Donagi result in \cite{Do} saying that if we have two polynomials $f, g \in S_d$, then $E'(f)=E'(g)$ implies that the corresponding hypersurfaces $V(f):f=0$ and $V(g):g=0$ are projectively equivalent. In fact, the general linear group $G=GL_{n+1}(\C)$ acts in an obvious way by substitutions on the vector space $S_d$ of homogeneous polynomials of degree $d$ and the tangent space at $f$ to the orbit $G \cdot f$ is exactly the vector space $J_{f,d}$, i.e.
\begin{equation}\label{tangent}
T_f(G\cdot f)=E'(f).
\end{equation}

In this last section we prove the following result, showing that in fact the tangent space $T_f(G\cdot f)$ determines the polynomial $f$ (up to a multiplicative constant factor), i.e. determines the hypersurface $V(f)$, in most cases.

\begin{thm}\label{thm2}
Suppose given two homogeneous polynomials $f,g\in \mathbb{C}[x_0,\cdots,x_n]$ such that $V(f) \ne V(g)$, i.e. $ f\notin\mathbb{C}^*g$, of degree $d>3$. Assume that the minimal degree of a syzygy involving the partial derivatives of $g$, denoted by $mdr_0(g)$, is at least  $3$.
Then, if  $T_f(G\cdot f)=T_g(G\cdot g)$, then either the hypersurface $V(g)\subset\mathbb{P}^n$ has a singularity of multiplicity $d-1$, or $g$ is of Sebastiani-Thom type (ST type).
\end{thm}

\begin{cor}
If  $T_f(G\cdot f)=T_g(G\cdot g)$ and $V(f) \ne V(g)$, with $V(g)$  smooth of degree $d>3$, then $g$ is of ST type.
\end{cor}

The definition of the invariant $mdr_0(g)$ is given in the last section.

The results of this note may have applications in extending the Torelli property in the sense of Dolgachev-Kapranov to singular hypersurfaces, see K. Ueda and M. Yoshinaga  \cite{UY} for the case of smooth hypersurfaces and  A. Dimca and E. Sernesi \cite{DS} for the case of some singular plane curves.

\bigskip

The author would like to thank A. Dimca for suggesting this problem and to Laboratoire J.-A. Dieudonn\'e for financial support and a stimulating working atmosphere.

\section{Preliminaries}

In this article, we  discuss homogeneous polynomials of a special type, namely  of Sebastiani-Thom type.

\begin{definition} \label{defST}
Given a homogeneous polynomial $f$ of degree $d$, we say $f$ is of Sebastiani-Thom type (ST type) if $f$ can be represented as
\[
f(x_0,\cdots,x_n)=g(x_0,\cdots,x_l)+h(x_{l+1},\cdots,x_n)
\]
for a choice of  homogeneous coordinates $(x_i)_{i=0}^n$ of $\mathbb{P}^n$, some $0\leq l\leq n-1$ and some degree $d$ homogeneous polynomials $g$ and $h$. When $h=0$ we say that $V(f)$ or $f$ is a cone.
\end{definition}

\begin{ex}
Any linear form
 and any quadratic form is of ST type if $n \geq 1$. That is why we  assume $d\geq 3$ in the sequel.
\end{ex}

Given homogeneous coordinates $(x_i),(x_i')$ of $\mathbb{P}^n$, and a homogeneous polynomial $f(x)=f'(x')$ of degree $d$, we  define
\[
f_i:=\frac{\partial f}{\partial x_i},\ \ f_{ij}:=\frac{\partial^2f}{\partial x_i\partial x_j},\ \ f_i':=\frac{\partial f'}{\partial x_i'};
\]
Note that $f_i$ is a polynomial in the $(x_i)$'s and $f'_i$ is a polynomial in the $(x'_i)$'s.

\begin{ex} \label{ex2}
If $f_i=0$ for some $i$, then $f$ is a cone, and hence of ST type. More generally, if $f_0,\cdots,f_n$ are linearly dependent, then $f$ is a cone.

We introduce some useful notations while we prove this easy fact. Consider the linear coordinate transformations
\[
x_i'=\sum_j a_{ij} x_j,\ \ x_i=\sum_jb_{ij}x_j',
\]
or
\[
X'=AX,\ \ X=BX',
\]
where $A=(a_{ij}), B=(b_{ij})$ are invertible matrices. Then under this transformation, we have
\[
f_j=\sum_i f_i'a_{ij},\ \ f_k'=\sum_lf_l b_{lk},
\]
that is,
\[
\begin{pmatrix}
f_0\\
\vdots\\
f_n
\end{pmatrix}
=A^t
\begin{pmatrix}
f_0'\\
\vdots\\
f_n'
\end{pmatrix}
\]
and
\[
\begin{pmatrix}
f_0'\\
\vdots\\
f_n'
\end{pmatrix}
=B^t
\begin{pmatrix}
f_0\\
\vdots\\
f_n
\end{pmatrix},
\]
where $A^t$ and $B^t$ denote the transposed matrices of $A$ and $B$ respectively.
We may formulate these two equations in a simpler way:
\[
\nabla f= A^t\nabla'f,\ \ \nabla'f=B^t\nabla f,
\]
where
\[
\nabla f=
\begin{pmatrix}
f_0\\
\vdots\\
f_n
\end{pmatrix}
\]
and
\[
\nabla' f=
\begin{pmatrix}
f_0'\\
\vdots\\
f_n'
\end{pmatrix}.
\]

Now if $f_0,\cdots, f_n$ are linearly dependent, then there exists an index $i$ such that
\[
f_i=a_0 f_0+\cdots+a_{i-1}f_{i-1}+a_{i+1}f_{i+1}+\cdots+a_nf_n.
\]
Let $A=(a_{ij})$ be defined as follows:
\[
A_{kl}=\delta_{kl},\ for \ k\neq i;
\]
and
\[
a_{ij}=-a_i,\ for\  j\neq i;\ \ a_{ii}=1.
\]
Then $A$ is invertible and
\[
A
\begin{pmatrix}
f_0\\
\vdots\\
f_n
\end{pmatrix}
=
\begin{pmatrix}
f_0\\
\vdots\\
f_{i-1}\\
0\\
f_{i+1}\\
\vdots\\
f_n
\end{pmatrix}.
\]
Let $X=BX'$ such that $B^t=A$. Then $f_i'=0$ and thus $f'$ does not depend on the variable $x_i'$, and hence $f$ is a cone.
\end{ex}

From the example above, we know that we can assume $f_0,\cdots,f_n$ are linearly independent when we discuss whether $f$ is of ST type or not.

Therefore, in the sequel we assume the following.

{\bf H1:} All given polynomials $f, g,\cdots$ are homogeneous of degree $d\geq3$;

{\bf H2:} The first order partial derivatives of the given polynomial $f$ are linearly independent,i.e. $f$ is not a cone.

\section{The general setting}

Given two polynomials $f$ and $g$, the following facts are obvious.

\noindent{\bf Fact1:} $J_f=J_g$ if and only if $E(f)=E(g)$.

\noindent{\bf Fact2:} If $V(f)$ is smooth, then $f_0,\cdots,f_n$ are linearly independent (see Example \ref{ex2})

\noindent{\bf Fact3:} Given $f$ and $g$ with $J_f=J_g$, then the singular schemes $Sing(V(f))$ and $Sing(V(g))$ are identical, as they are defined by the ideal $J_f=J_g$. In particular,
if $V(f)$ is smooth, then $V(g)$ is smooth.

\bigskip

Now we consider two polynomials $f$ and $g$ with the condition that $J_f=J_g$ and $f_0,\cdots,f_n$ are linearly independent. Then $g_0,\cdots,g_n$ are also linearly independent and there exists an invertible matrix $C$ such that
\[
\begin{pmatrix}
f_0\\
\vdots\\
f_n
\end{pmatrix}
=C
\begin{pmatrix}
g_0\\
\vdots\\
g_n
\end{pmatrix}
\]
or
\[
\nabla f=C\nabla g.
\]

From linear algebra, we know that there is an invertible matrix $P$ and a Jordan matrix $J$ such that
\[
C=P^{-1}JP
\]
where
\[
J=
\begin{pmatrix}
J_1 & &  & \\
   & J_2 & &\\
   &  & \ddots &\\
   &  &     & J_k
\end{pmatrix}
\]
with
\[
J_i=
\begin{pmatrix}
\lambda_i & 1 & 0 & \cdots & \cdots & 0\\
          & \lambda_i & 1 & 0 &\cdots  &0 \\
          & & \ddots &\ddots  &\ddots  &\vdots\\
          & & & \lambda_i& 1 &0\\
          & & & & \lambda_i &1\\
          & & & & &\lambda_i
\end{pmatrix}.
\]
So we have
\[
\nabla f=P^{-1}JP\nabla g,
\]
that is
\[
P\nabla f=JP\nabla g,
\]
i.e.
\[
\nabla'f=J\nabla'g,
\]
where $X=P^tX'$ is the corresponding linear coordinate transformation.

Thus without loss of generality (WLOG), we can assume the following.

\noindent{\bf H3:} The given polynomials $f$ and $g$ satisfy $\nabla f=J\nabla g$.

\noindent{\bf H4:} $V(f)\neq V(g)$.

The last condition is equivalent to the condition that the matrix $J$ is not of the form $\lambda I_{n+1}$.

\section{Two key cases}
\subsection{$J$ is diagonal, i.e. $k=n+1$}
If $J=diag\{\lambda_0,\cdots,\lambda_n\}$, then
$
f_i=\lambda_i g_i,
$
so
\[
f_{ij}=(f_i)_j=\lambda_i g_{ij}
\text{ and }
f_{ij}=f_{ji}=(f_j)_i=\lambda_j g_{ji}=\lambda_j g_{ij}.
\]
So we get
\[
(\lambda_i-\lambda_j)g_{ij}=0,
\]
implying that
\[
g_{ij}=0\ \ whenever\ \lambda_i\neq\lambda_j.
\]
After some coordinate renumbering, we can rewrite $J$ as
\[
J=diag\{\mu_1 I_{n_1},\mu_2 I_{n_2},\cdots, \mu_pI_{n_p}\},
\]
where $\mu_1,\cdots,\mu_p$ are distinct. The above implies that $g$ can be represented as
\[
g(x_0,\cdots, x_n)=H(x_0,\cdots, x_{n_1-1})+K(x_{n_1},\cdots,x_{n_1+n_2-1})+\cdots+L(x_{n+1-n_p},\cdots,x_n).
\]
Since $J$ is not of the form $\lambda I_{n+1}$, $1\leq n_1<n+1$, so $g$ (and $f$ by symmetry) is of ST type.

\subsection{$J$ is a Jordan block, i.e. $k=1$}

If $k=1$, then $\nabla f=J\nabla g$ implies that for some $\lambda\in\mathbb{C^*}$ one has the the following.
\begin{equation}
\label{f,g}
f_0=\lambda g_0+g_1,\cdots, f_{n-1}=\lambda g_{n-1}+g_n,\ \ f_n=\lambda g_n.
\end{equation}

For $0\leq r,s\leq n-1$, we have
\[
f_r=\lambda g_r+g_{r+1},
\text{ and }
f_s=\lambda g_s+g_{s+1},
\]
so we get
\[
f_{rs}=(f_r)_s=\lambda g_{rs}+g_{r+1 s},
\text{ and }
f_{sr}=(f_s)_r=\lambda g_{sr}+g_{s+1 r}.
\]
Hence $f_{rs}=f_{sr}, g_{rs}=g_{sr}$ implies that
\[
g_{r+1 s}=g_{s+1 r}.
\]
Moreover,
\[
f_{rn}=(f_r)_n=\lambda g_{rn}+g_{r+1 n}
\text{ and }
f_{nr}=(f_n)_r=\lambda g_{nr},
\]
imply
\[
g_{r+1 n}=0,
\]
that is $g_{1n}=g_{2n}=\cdots=g_{nn}=0$, so $g$ is of the form
\[
g(x_0,\cdots, x_n)= ax_0^{d-1}x_n+h(x_0,\cdots,x_{n-1}),
\]
for some $a \in \C^*$.
Note that in this case $V(g)$ must be singular, since $p=[0:\cdots:0:1]$ is a singular point of $V(g)$ with multiplicity $d-1$.

Since $g_{r+1, s}=g_{s+1, r}=g_{r, s+1}$, we have
\[
0=g_{1, n}=g_{2 ,n-1}=\cdots=g_{n-1, 2}
\]
\[
0=g_{2 ,n}=g_{3 ,n-1}=\cdots=g_{n-1, 3},
\]
\[
0=g_{n-1, n}=g_{n, n-1},
\]
\[
g_{nn}=0,
\]
So
$
g_{ij}=0$ for $ i+j\geq n+1$ and
hence the Hessian matrix of $g$ is of the type
\begin{equation}
\label{Hess}
\begin{pmatrix}
* & * & * & * & g_{0n}\\
* & * & * & g_{0n} & 0\\
* & * &\iddots& \iddots & \vdots\\
* &g_{0n}&0 &\cdots& 0\\
g_{0n} & 0 &\cdots &\cdots & 0
\end{pmatrix}.
\end{equation}

One has
\[
g_{0n}=(d-1)ax_0^{d-2},
\]
so the determinant of Hessian matrix is
\[
Hess(g)=c x_0^{(d-2)(n+1)},
\]
for some constant $c \in \C^*$.
We also know that under linear coordinate transformations, the determinant of the Hessian matrix changes by multiplication by a nonzero constant. Suppose  in some linear coordinates $(x_i')$, $g$ can be written as
\[
g(x_0',\cdots,x_n')=h(x_0',\cdots,x_l')+k(x_{l+1}',\cdots,x_n').
\]
Then $Hess(g)=Hess(h)\times Hess(k)$. But $Hess(h)$ is a polynomial $P(x_0',\cdots,x_l')$ and $Hess(h)$ is a polynomial $Q(x_{l+1}',\cdots,x_n')$; both $P$ and $Q$ are of strictly positive degrees. Since $Hess(g)$ is a power of a linear form in the $x_i'$, this is a contradiction, and hence
 $g$ is not of ST type.

\begin{ex} \label{ex} First we treat the case $n=2$, $d=3$.
In the homogeneous coordinates $(x,y,z)$, set
\[
g(x,y,z)=\frac{1}{2}(x^2 z+x y^2),
\]
then
\[
\nabla g=
\begin{pmatrix}
xz+\frac{1}{2}y^2\\
xy\\
\frac{1}{2} x^2
\end{pmatrix}
\]
then the hessian matrix of $g$ is
\[
\begin{pmatrix}
z & y & x\\
y & x & 0\\
x & 0 & 0\\
\end{pmatrix}
\]
i.e. of the form given by equation \eqref{Hess}.

Now let
\[
f(x,y,z)=\frac{1}{2}(x^2 z+ xy^2+x^2 y),
\]
then
\[
\nabla f=
\begin{pmatrix}
xz+\frac{1}{2} y^2+xy\\
xy+\frac{1}{2}x^2\\
\frac{1}{2}x^2\\
\end{pmatrix}
=
\begin{pmatrix}
1 & 1 & 0\\
0 & 1 & 1\\
0 & 0 & 1
\end{pmatrix}
\nabla g,
\]
exactly as in the equations \eqref{f,g}. In fact, we determine first the $g_i$ using these equations, and then we obtain $f$ by integration. It follows that $\mathbb{C}\langle f_x,f_y,f_z\rangle=\mathbb{C}\langle g_x,g_y,g_z\rangle$,  $V(f)\neq V(g)$, and $g$ is not of the ST type. Symmetrically, $f$ is not of ST type. The point $p=[0:0:1]$ has multiplicity $d-1=2$.
Similarly, we can also consider $f_\lambda(x,y,z)=\frac{1}{2}(\lambda x^2 z+\lambda xy^2+x^2 y)$ with $\lambda\neq0$.
\end{ex}

\begin{ex}\label{existence}
It is not easy to explicitly construct a $g\in S_{n,d}$ which is not of ST type for general $d\geq3$ and $n$ though, of great importance is the existence of such functions. In Corollary \ref{corA}, we use this to show that a generic polynomial is not of ST type.

Here is a simple proof of the existence. As usual, we see $S_{n,d}$ as a vector space in the natural way. For $g\in S_{n,d}$, we have proved that if $J\nabla g=\nabla f$ for some $f\in S_{n,d}$, then $g$ is not of ST type. Note that here $J$ is a $(n+1)\times(n+1)$ Jordan block with diagonal entries equal to $1$.

We construct such functions inductively on $n$; more precisely, if there exists $h\in S_{n-1,d}=\mathbb{C}[x_0,\cdots,x_{n-1}]_d$ not of ST type as above with
$$h_{r+1s}=h_{r,s+1}, h_{r+1,n-1}=0,\ 0\leq r,s<n-1,$$
implying that $h$ is of the form $x_0^{d-1}x_{n-1}+K(x_0,\cdots,x_{n-2})$, then we can construct $g\in S_{n,d}$ by setting
\begin{equation}
\label{general}
g=\int_0^{x_1}h_0 +\int_0^{x_1} h_1|_{x_1=0}+\cdots+\int_0^{x_{i+1}} h_i|_{x_1=\cdots=x_i=0}+\int_0^{x_n}h_{n-1}|_{x_1=\cdots=x_{n-1}=0},
\end{equation}
where $h_i|_{x_1=\cdots=x_i=0}$ denotes the polynomial obtained from $h_i$ by setting $x_1=\cdots=x_i=0$ and $\int_0^{x_{i+1}}$ denotes the definite integral with respect to the $(i+1)$-th variable from $0$ to $x_{i+1}$, meaning that
\[
\int_0^{x_{i+1}}h_i|_{x_1=\cdots=x_i=0}:=\int_0^{x_{i+1}}h_i(x_0,0,\cdots,0, y_{i+1}, x_{i+2},\cdots,x_{n-1})\,\mathrm{d}y_{i+1}.
\]
We can check that $g$ given by \eqref{general} satisfies $g_{r+1}=h_r$ for $0\leq r\leq n-1$, and in fact, we get the formula by integration from these conditions. The corresponding $f$ is $g+h$.

If $n=2$ and $d\geq3$, then exactly as $n=2,d=3$ we get
\begin{equation}
\label{n=2}
g=x^{d-1}z+\frac{d-1}{2}x^{d-2}y^2 \text{ and } f=g+x^{d-1}y.
\end{equation}
Using $\eqref{general}$, if $n=3$, then we get in the coordinates $x,y,z,t$
\begin{equation}
\label{n=3}
g=x^{d-1}t+(d-1)x^{d-2}yz +\frac{(d-1)(d-2)}{6}x^{d-3}y^3 \text{ and }  f=g+x^{d-1}z+\frac{d-1}{2}x^{d-2}y^2.
\end{equation}
Similarly, if $n=4$, then we get in the coordinates $x,y,z,t,u$
$$
g=x^{d-1}u+(d-1)x^{d-2}(yt +1/2z^2)+\frac{(d-1)(d-2)}{2}x^{d-3}y^2z+\frac{(d-1)(d-2)(d-3)}{24}x^{d-4}y^4$$
 and
$$ f=g+x^{d-1}t+(d-1)x^{d-2}yz+\frac{(d-1)(d-2)}{6}x^{d-3}y^3.$$
\end{ex}

\begin{ex} Now we consider special case $n=2$ for any $d$ with coordinates $x,y,z$. We start from the equation $\nabla f=J\nabla g$, where $J=diag\{J_1,\cdots,J_k\}$ is a Jordan matrix with only one eigenvalue 1 and $J_i$ is a Jordan block. We will determine all possibilities of $f$ and $g$. And we will discuss 2 cases: $k=1$ and $k=2$.

{\bf Case 1: $k=1.$} We have
\[
J=
\begin{pmatrix}
1 & 1 & 0\\
0 & 1 & 1\\
0 & 0 & 1
\end{pmatrix}
\]
Then just as what we have discussed, we get
\[
g_{yz}=g_{zz}=0,\ g_{yy}=g_{xz}.
\]
From $g_{yz}=g_{zz}=0$, we have
\[
g(x,y,z)=x^{d-1}z+K(x,y),
\]
and from $g_{yy}=g_{xz}$, we have $K_{yy}=(d-1)x^{d-2}$, so by integration, we get
\[
K(x,y)=\frac{d-1}{2}x^{d-2}y^2+A x^{d-1}y+B x^d,
\]
thus
\[
g(x,y,z)=x^{d-1} z+\frac{d-1}{2}x^{d-2}y^2+A x^{d-1} y+B x^d,
\]
and
\[
f(x,y,z)=g(x,y,z)+x^{d-1}y+\frac{A}{d} x^d,
\]
where $A,B$ are constants which can be arbitrarily chosen. By setting $A=B=0$, we get the results in Example \ref{ex} for $n=2$ and $d\geq3$.

It is easy to see that in this case $x^{d-2}|g$ and $x^{d-2}|f$.

{\bf Case 2: $k=2$.} We have
\[
J=
\begin{pmatrix}
1 & 1 & 0\\
0 & 1 & 0\\
0 & 0 & 1
\end{pmatrix}
\]
so
\[
f_x=g_x+g_y,\ \ f_y=g_y,\ \ f_z=g_z,
\]
from $f_{xy}=f_{yx}, f_{xz}=f_{zx}$ we get
\[
g_{yy}=0,\ \ g_{yz}=0,
\]
so
\[
g(x,y,z)=x^{d-1} y+H(x,z),
\]
and
\[
f(x,y,z)=g(x,y,z)+\frac{1}{d}x^d.
\]
Here $H(x,z)\in\mathbb{C}[x,z]_d$ can be arbitrarily chosen.

 This case is very different from the first one. First, we do not have the divisibility as in case 1. Also, as we have proved, any $g$ in case 1 is not of ST type though, we can have $g$ of ST type in this case, for instance, by setting $H(x,z)=z^d$ and get $g'=x^{d-1} y+z^d$.
In fact one can show that $g$ is of ST type if and only if $g$ is projectively equivalent to $g'$ above.

\end{ex}

\section{The general case}

We start from the equation $\nabla f=J\nabla g$, and consider the decomposition $J=diag\{J_1,\cdots,J_k)$ where $J_i$ is a $n_i\times n_i$ Jordan block.

If all $n_i=1$, then $J$ is diagonal and this case was discussed above.

If $J$ is not diagonal, then there exists $n_i>1$ and WLOG, we assume $n_1>1$. Set $N_1=n_1-1$ and note that exactly as above we get
\[
f_{N_1-1}=\lambda_1 g_{N_1-1}+g_{N_1}
\text{ and }
f_{N_1}=\lambda_1 g_{N_1}.
\]
Then from $f_{N_1-1 N_1}=f_{N_1 N_1-1}$, we have $g_{N_1 N_1}=0$, thus
\[
g(x_0,\cdots,x_n)=x_{N_1}H(x_0,\cdots,\widehat{x_{N_1}},\cdots,x_n)+K(x_0,\cdots,\widehat{x_{N_1}},\cdots,x_n),
\]
implying that $V(g)$ has a singularity of multiplicity $d-1$ at the point having $x_{N_1}=1$ and all the other coordinates equal to zero.

\section{An application}

We set  $S_{n,d}=\mathbb{C}[x_0,\cdots,x_n]_d$, the space of homogeneous polynomials of degree $d\geq 3$ in $(n+1)$ variables, and $\mathbb{P}(S_{n,d})$  the corresponding projective space. Let $\mathscr{G}$ be the subset of $g\in \mathbb{P}(S_{n,d})$ such that either $g$ is  of ST type or $V(g)$ has some singularities of multiplicity $d-1$; let $\mathscr{U}=\mathbb{P}(S_{n,d})\setminus\mathscr{G}$.

We show below that $\mathscr{G}$ is a proper closed algebraic subset of $\mathbb{P}(S_{n,d})$, or equivalently, $\mathscr{U}$ is a nonempty open subset of $\mathbb{P}(S_{n,d})$ in the Zariski topology. Note that by Theorem \ref{thm1}, any $g\in\mathscr{U}$ has the following property: if $f\in\mathbb{P}(S_{n,d})$ is a polynomial whose  first order derivatives $f_0=\frac{\partial f}{\partial x_0},\cdots,f_n=\frac{\partial f}{\partial x_n}$ form a basis of the vector space $E(g)$, then $f=g$ in $\mathbb{P}(S_{n,d})$. More precisely, we have the following result.

\begin{cor} \label{corA}
The set  $\mathscr{U}$ is a nonempty Zariski open subset  of $\mathbb{P}(S_{n,d}),$ when $ d\geq 3$.  Moreover, for any $g\in\mathscr{U}$, there exists a unique basis $F_0,\cdots,F_{n}$ for the vector space $E(g)$, up to a common nonzero multiple, satisfying
\[
\frac{\partial F_i}{\partial x_j}=\frac{\partial F_j}{\partial x_i}
\]
for all $0 \leq i <j \leq n$.
\end{cor}

This is just the celebrated result due to J. Carlson and Ph. Griffith, see section 4, (b) in \cite{CG}, but without the explicit description of $\mathscr{U}$ given in our Corollary. In \cite{CG}, this result is proved by an implicit use of the  ``Curve Selection Lemma", see for instance \cite{Mi} for a statement and proof of this Lemma. As a result, one cannot get  information about the precise meaning of ``generic"  from their proof. In this paper, they also conjecture that the map $g\mapsto E(g)$  from $\mathbb{P}(S_{n,d})$ to an obvious Grassmannian is injective as long as $V(g)$ does not possess singularities which are ``large in relation to the degree". And indeed, by Theorem \ref{thm1}, we have confirmed this conjecture.

\subsubsection{Proof of Corollary \ref{corA}}
Now we prove that $\mathscr{G}$ is a proper closed algebraic subset of $\mathbb{P}(S_{n,d})$. Set $\mathscr{G}_1$ be the subset of $g\in\mathbb{P}(S_{n,d})$ which is of ST type, and  $\mathscr{G}_2$ be the subset of $g\in\mathbb{P}(S_{n,d})$ with $V(g)$ admitting singularities of multiplicity $\geq d-1$. Since a cone $f$ is of ST type, we know that $\mathscr{G}=\mathscr{G}_1\cup\mathscr{G}_2$, thus it suffices to prove that $\mathscr{G}_1$ and $\mathscr{G}_2$ are proper closed algebraic subsets in $\mathbb{P}(S_{n,d})$.

First we know that $g\in S_{n,d}$ is of ST type if and only if there exists an integer $l$, an invertible matrix $A$ and a pair $(h,k)\in S_{l,d}\times S_{n-l-1,d}$ such that $0\leq l<n$ and
\begin{equation}
\label{eqst1}
g(x_0,\cdots,x_n)=h(x_0',\cdots,x_l')+k(x_{l+1}',\cdots,x_n'),
\end{equation}
where $X'=AX$. This leads us to set
\[
\mathscr{H}_1(l):=\{(g,(h,k),A)\in\mathbb{P}(S_{n,d})\times\mathbb{P}(S_{l,d}\times S_{n-l-1,d})\times PGL_{n+1}: \eqref{eqst1}\ \text{holds}\},
\]
for $0\leq l<n$. Note that in this definition, $PGL_{n+1}$ is the set of projective transformations of $\mathbb{P}^n$, and $A\in PGL_{n+1}$ naturally induces a linear coordinates transformation $X'=AX$ under which \eqref{eqst1} holds.

Let $Mat_{n+1,n+1}$ be the vector space of all $(n+1)\times (n+1)$ matrices. We see $\mathscr{H}_1(l)$ as the subset of $\mathbb{P}(S_{n,d})\times\mathbb{P}(S_{l,d}\times S_{n-l-1,d})\times\mathbb{P}(Mat_{n+1,n+1})$, then $\mathscr{H}_1(l)$ is constructible.

Let
\[
\pi_l:\mathbb{P}(S_{n,d})\times\mathbb{P}(S_{l,d}\times S_{n-l-1,d})\times\mathbb{P}(Mat_{n+1,n+1})\rightarrow\mathbb{P}(S_{n,d})
\]
be the projection onto the first factor, then $\pi_l$ is a morphism; by Chevalley's theorem, $\pi_l(\mathscr{H}_1(l))$ is a constructible subset of $\mathbb{P}(S_{n,d})$. Thus
\[
\mathscr{G}_1=\bigcup_{l=0}^{n-1}\pi_l(\mathscr{H}_1(l))
\]
is also constructible.

To prove that $\mathscr{G}_1$ is a closed algebraic subset of $\mathbb{P}(S_{n,d})$, it suffices to show that it is a closed subset in the classical topology. To see this, let $\{g_i\}_{i=1}^\infty\in\mathscr{G}_1$ and $g_i\rightarrow g_\infty\in\mathbb{P}(S_{n,d})$, we should prove that $g_\infty\in\mathscr{G}_1$, i.e. $g_\infty$ is of ST type.

Since $g_i$ is of ST type, by definition, there exists an integer $l_i$, a pair $(h_i,k_i)\in S_{l_i,d}\times S_{n-l_i-1,d}$ and an invertible matrix $A_i$ inducing a coordinate transformation $X'^{(i)}=A_iX$ such that
\begin{equation}
\label{eqst2}
0\leq l_i<n\ \text{and}\ g_i(x_0,\cdots,x_n)=h_i(x_0'^{(i)},\cdots,x_{l_i}'^{(i)})+k_i(x_{l_i+1}'^{(i)},\cdots,x_n'^{(i)}).
\end{equation}
WLOG, we may assume $l_i=l$ with $0\leq l<n$ for $i\geq1$.

We see $((h_i,k_i),A_i)$ as elements of $\mathbb{P}(S_{l,d}\times S_{n-l-1,d})\times\mathbb{P}(Mat_{n+1,n+1})$, then by the compactness of the latter, we may assume this sequence converges, so there exists $((h_\infty,k_\infty),A_\infty)\in\mathbb{P}(S_{l,d}\times S_{n-l-1,d})\times \mathbb{P}(Mat_{n+1,n+1})$ such that
\[
((h_i,k_i),A_i)\rightarrow ((h_\infty,k_\infty),A_\infty)\ \text{in}\ \mathbb{P}(S_{l,d}\times S_{n-l-1,d})\times\mathbb{P}(Mat_{n+1,n+1});
\]
since $g_i\rightarrow g_\infty$ in $\mathbb{P}(S_{n,d})$, letting $i\rightarrow\infty$ in \eqref{eqst2}, we have
\begin{equation}
\label{eqst3}
g_\infty(x_0,\cdots,x_n)=h_\infty(x_0',\cdots,x_l')+k_\infty(x_{l+1}',\cdots,x_n'),\ \text{in}\ \mathbb{P}(S_{n,d})
\end{equation}
with $X'=A_\infty X$.

If $A_\infty$ is invertible, we see that $g_\infty$ is of ST type. Otherwise, $A_\infty$ is singular, then from \eqref{eqst3}, the first order derivatives of $g_\infty$ are linearly dependent, so $g_\infty$ is a cone. From Example \ref{ex2}, we know that $g_\infty$ is of ST type, i.e. $g_\infty\in\mathscr{G}_1$. Thus, $\mathscr{G}_1$ is a closed subset in the classical topology, implying that it is a closed algebraic subset of $\mathbb{P}(S_{n,d})$. From Example \ref{existence}, we know that $\mathscr{G}_1$ is a proper subset.

As the final step of the proof, we  show that $\mathscr{G}_2$ is a proper closed algebraic subset of $\mathbb{P}(S_{n,d})$. Set
\[
\mathscr{H}_2=\{(g,x)\in\mathbb{P}(S_{n,d})\times\mathbb{P}^n: x\ \text{is a singular point of $V(g)$ with multiplicity $\geq d-1$}\},
\]
then $\mathscr{H}_2$ is a closed algebraic subset in $\mathbb{P}(S_{n,d})\times\mathbb{P}^n$. Let
\[
pr: \mathbb{P}(S_{n,d})\times\mathbb{P}^n\rightarrow\mathbb{P}(S_{n,d})
\]
be the projection onto the first factor, then $\mathscr{G}_2=pr(\mathscr{H}_2)$ is closed and algebraic, since $pr$ is closed. $\mathscr{G}_2$ does not contain smooth $g's$, so $\mathscr{G}_2$ is also a proper subset.

\subsubsection{Remark}
 We can also show that $\mathscr{G}_1$ is a proper subset in another way. Note that if a smooth hypersurface $V(f)$ is defined by a polynomial $f$ of ST type as in Definition \ref{defST}, then the linear automorphism group $Lin(V(f))$ of the hypersurface $V(f)$ contains a cyclic subgroup of order $d$, consisting of the elements $u_a: X \mapsto X'$ with $x'_i=x_i$ for $i \leq l$ and $x'_i=ax_i$ for $i>l$ and $a^d=1$. It is known that this linear automorphism group $Lin(V(f))$ of a  hypersurface is trivial for $n \ge 3$ and $d \ge 3$ when the coefficients of $f$ are algebraically independent over  $\mathbb{Q}$ , see Theorem 5 in \cite{MM}. Hence such a hypersurface is not of ST type. The remaining  case $n=2$ follows from Example \ref{existence}, where we have constructed homogeneous polynomials not of ST type for this value of $n$ and any degree $d \geq 3$.

\section{Hypersurfaces determined by tangent spaces}

Let $f$ be a homogeneous polynomial of degree $d$ in the polynomial ring $S$ and denote by $f_0,...,f_n$ the corresponding partial derivatives.

One can consider the graded $S-$submodule $AR(f) \subset S^{n+1}$ of {\it all relations} involving the $f_j$'s, namely
$$a=(a_0,...,a_n) \in AR(f)_m$$
if and only if  $a_0f_0+a_1f_1+...+a_nf_n=0$.

\begin{definition} \label{mdr}
We define the minimal degree of a syzygy involving the partial derivatives of $f$ to be the positive integer
$$mdr_0(f) = \min \{m \  \ : \  \ AR(f)_m \ne 0 \} .$$

\end{definition}

\begin{ex} \label{mdrex1}

If $V(f)$ is a smooth hypersurface of degree $d$, then $f_0,...,f_n$ is a regular sequence in $S$ and hence $mdr_0(f)=d-1.$

\end{ex}

\begin{ex} \label{mdrex1}

If $V(f)$ is a nodal plane curve of degree $d$, then  $mdr_0(f)=d-2$ if the curve is reducible and
$mdr_0(f)=d-1$ if the curve is irreducible, see \cite{DSt1}, \cite{DSt2}. In particular our results above apply to nodal curves of degree at least $5$ (resp. to irreducible nodal curves of degree at least $4$).

If $V(f)$ is a plane curve having only nodes and cusps, then $mdr_0(f) \geq 5d/6-2$ by \cite{DS}, hence our results apply to such curves if their degree $d$ is at least $6$.
\end{ex}

\begin{ex} \label{mdrex2}

If $V(f)$ is a nodal surface of degree $d$ in $\PP^3$, then
$mdr_0(f)=d-1$  \cite{D1}, \cite{DSa}. In particular our results above apply to nodal surfaces of degree at least $4$.

If $V(f)$ is a nodal threefold in $\PP^4$, then
$mdr_0(f)=d-1$ by \cite{DSt1}, hence our results apply to such threefolds if their degree $d$ is at least $4$.
\end{ex}

\subsection{The proof of   Theorem \ref{thm2}}

Now we assume that we are given two homogeneous polynomials $f,g\in S_d$ with $E'(f)=E'(g)$. Then for some $a_{ij,\alpha\beta}\in\mathbb{C}$,
\[
x_if_j=\sum_{\alpha,\beta} a_{ij,\alpha\beta}x_\alpha g_\beta.
\]
so
\[
x_ix_jf_k=x_i\sum_{\alpha,\beta}a_{jk,\alpha\beta}x_\alpha g_\beta=\sum_{\alpha,\beta}a_{jk,\alpha\beta}x_ix_\alpha g_\beta,
\]
\[
x_ix_jf_k=x_j(x_if_k)=x_j\sum_{\alpha,\beta}a_{ik,\alpha\beta}x_\alpha g_\beta=\sum_{\alpha,\beta}a_{ik,\alpha\beta}x_jx_\alpha g_\beta.
\]
Thus,
\[
\sum_{\alpha,\beta}a_{jk,\alpha\beta}x_ix_\alpha g_\beta=\sum_{\alpha,\beta}a_{ik,\alpha\beta}x_jx_\alpha g_\beta
\]
that is,
\[
\sum_{\beta}\bigg(\sum_{\alpha}(a_{jk,\alpha\beta}x_i-a_{ik,\alpha\beta}x_j)x_\alpha\bigg)g_\beta=0.
\]
Since by assumption, $mdr_0(g) \geq 3$,  implying that for all $i,j,k,\beta$
\[
\sum_{\alpha}(a_{jk,\alpha\beta}x_i-a_{ik,\alpha\beta}x_j)x_\alpha=0.
\]
Since $x_0,\cdots,x_n$ is a regular sequence, so for all $i,j,k,\alpha,\beta$,
\[
a_{jk,\alpha\beta}x_i-a_{ik,\alpha\beta}x_j\in(x_0,\cdots,\widehat{x_\alpha},\cdots,x_n).
\]
Now for any given $\beta,k$, we fix $\alpha$; let $i=\alpha, j\neq \alpha$, we get
\[
a_{jk,\alpha\beta}x_\alpha-a_{\alpha k,\alpha\beta}x_j\in(x_0,\cdots,\widehat{x_\alpha},\cdots,x_n),
\]
so $a_{jk,\alpha\beta}=0$. Therefore, $a_{jk,\alpha\beta}=0$ for $j\neq\alpha$.

Since as above, $a_{jk,\alpha\beta}=0$ for $j\neq\alpha$,
\[
x_jf_k=\sum_{\alpha,\beta}a_{jk,\alpha\beta}x_\alpha g_\beta=\sum_\beta a_{jk,j\beta}x_jg_\beta,
\]
thus,
\[
f_k=\sum_\beta a_{jk,j\beta}g_\beta.
\]
So $E(f)=E(g)$, and the proof is completed by using the main result in \ref{thm1}.

\subsection{An application}
Let $g \in \PP(S_d)$ be a polynomial such that $V(g)$ is smooth of degree $d >3$  and $g$ is not of Sebastiani-Thom type. By the proof of Corollary \ref{corA}, such $g$'s give a nonempty Zariski open subset $\mathscr{U}'$ of $\mathbb{P}(S_d)$.

Consider the map $\varphi: f\mapsto E'(f)$  from $\mathscr{U}'$ to an obvious Grassmannian.
Then Theorem \ref{thm2} implies the following.

\begin{cor}\label{corB}

With the above assumptions, the map  $\varphi$ is injective.

\end{cor}

Similar results apply if the condition of smoothness of $V(g)$ is  replaced by  the condition $mdr_0(g) \geq 3$ (as in Examples \ref{mdrex1}, \ref{mdrex2}) and $V(g)$ has no singularity of multiplicity $d-1$.

In conclusion, Corollary \ref{corA} and Corollary \ref{corB} imply that a generic $f\in\mathbb{P}(S_d)$ can be reconstructed from either of the two homogeneous components $J_{f,d-1}$ and $J_{f,d}$ of the Jacobian ideal $J_f$. 
One can also consider similar questions for other homogeneous components of $J_f$.

\end{document}